\documentclass[amstex,12pt,english,amssymb]{article}

\usepackage{mathtext}
\usepackage[cp1251]{inputenc}
\usepackage[T2A]{fontenc}
\usepackage[english]{babel}
\usepackage[dvips]{graphicx}
\usepackage{amsmath}
\usepackage{amssymb}
\usepackage{amsxtra}
\usepackage{latexsym}
\usepackage{ifthen}

\usepackage[dvips]{graphicx}
\usepackage{epstopdf}
\usepackage{epsfig}
\usepackage{epic}
\usepackage{eepic}

\begin{document}

\numberwithin{equation}{section}
\newcommand{\abs}[1]{\lvert#1\rvert}
\newcommand{\blankbox}[2]{%
  \parbox{\columnwidth}{\centering
    \setlength{\fboxsep}{0pt}%
    \fbox{\raisebox{0pt}[#2]{\hspace{#1}}}
} }

\newcounter{lemma}[section]
\newcommand{\lemma}{\par \refstepcounter{lemma}%
{\bf Lemma \arabic{section}.\arabic{lemma}.}}
\renewcommand{\thelemma}{\thesection.\arabic{lemma}}

\newcounter{corollary}[section]
\newcommand{\corollary}{\par \refstepcounter{corollary}%
{\bf Corollary \arabic{section}.\arabic{corollary}.}}
\renewcommand{\thecorollary}{\thesection.\arabic{corollary}}

\newcounter{remark}[section]
\newcommand{\remark}{\par \refstepcounter{remark}%
{\bf Remark \arabic{section}.\arabic{remark}.}}
\renewcommand{\theremark}{\thesection.\arabic{remark}}

\newcounter{theorem}[section]
\newcommand{\theorem}{\par \refstepcounter{theorem}%
{\bf Theorem \arabic{section}.\arabic{theorem}.}}
\renewcommand{\thetheorem}{\thesection.\arabic{theorem}}

\newcounter{proposition}[section]
\newcommand{\proposition}{\par \refstepcounter{proposition}%
{\bf Proposition \arabic{section}.\arabic{proposition}.}}
\renewcommand{\theproposition}{\thesection.\arabic{proposition}}

\renewcommand{\theequation}{\arabic{section}.\arabic{equation}}

\def\kohta #1 #2\par{\par\noindent\rlap{#1)}\hskip30pt
\hangindent30pt #2\par}
\def\A{{{\cal {A}}}}
\def\C{{{\Bbb C}}}
\def\R{{{\Bbb R}}}
\def\Rn{{{\Bbb R}^n}}
\def\lC{{\overline {{\Bbb C}}}}
\def\lRn{{\overline {{\Bbb R}^n}}}
\def\lRm{{\overline {{\Bbb R}^m}}}
\def\lRk{{\overline {{\Bbb R}^k}}}
\def\lBn{{\overline {{\Bbb B}^n}}}
\def\Bn{{{\Bbb B}^n}}
\def\B{{{\Bbb B}}}
\def\dist{{{\rm dist}}}
\let\text=\mbox

\def\Xint#1{\mathchoice
   {\XXint\displaystyle\textstyle{#1}}%
   {\XXint\textstyle\scriptstyle{#1}}%
   {\XXint\scriptstyle\scriptscriptstyle{#1}}%
   {\XXint\scriptscriptstyle\scriptscriptstyle{#1}}%
   \!\int}
\def\XXint#1#2#3{{\setbox0=\hbox{$#1{#2#3}{\int}$}
     \vcenter{\hbox{$#2#3$}}\kern-.5\wd0}}
\def\dashint{\Xint-}

\def\cc{\setcounter{equation}{0}
\setcounter{figure}{0}\setcounter{table}{0}}

\overfullrule=0pt

\title{{\bf On a variational method \\ for the Beltrami equations}}

\author{{\sc  Tatyana Lomako, Vladimir Ryazanov}\\}
\date{\today \hskip 4mm ({\tt LR16042011.tex})}
\maketitle

\begin{abstract}
We construct variations for the classes of regular solutions to
degenerate Beltrami equations with restrictions of the set-theoretic
type for the complex coefficient. On this basis, we prove the
variational maximum principle and other necessary conditions of
extremum.
\end{abstract}



\large

\section{Introduction}

Let $D$ be a domain in ${\Bbb C}$, $\overline{{\Bbb C}}\,= \,{{\Bbb
C}}\,\cup\,\left\{\infty\right\}$. The {\bf Beltrami equation} in
$D$ has the form
\begin{equation}\label{eq1}
f_{\overline{z}}=\mu(z)\cdot f_{z}
\end{equation}
where $\mu(z): D\rightarrow \mathbb{C}$ is a measurable function
with $|\mu(z)|<1$ a.e., $f_{\overline{z}}=\overline{\partial}
f=\left(f_{x}+if_{y}\right)/2$, $f_{z}=\partial
f=\left(f_{x}-if_{y}\right)/2$, $z=x+iy$, $f_{x}$ and $f_{y}$ denote
the partial derivatives of the mapping $f$ in $x$ and $y$,
respectively. The function $\mu$ is the {\bf complex coefficient}
and
\begin{equation}\label{eq2}
K_{\mu}(z)=\frac{1+|\mu(z)|}{1-|\mu(z)|}
\end{equation}
is the {\bf  dilatation quotient} or simply the {\bf  dilatation} of
equation (\ref{eq1}).

Recall that a mapping $f:D\rightarrow \mathbb{C}$ is called {\bf
regular at a point} $z_{0}\in D$ if $f$ has a total differential at
the point and its Jacobian
$J_{f}(z)=|f_{z}|^{2}-|f_{\overline{z}}|^{2}\neq0$ (see, e.g., I.1.6
in \cite{LV}). A homeomorphism $f$ of the class $W_{loc}^{1,\,1}$ is
called {\bf regular} if $J_{f}(z)>0$ a.e. Finally, the {\bf regular
solution} of the Beltrami equation (\ref{eq1}) in the domain $D$ is
a regular homeomorphism that satisfies (\ref{eq1}) a.e. in $D$. The
notion of the regular solution was first introduced in the paper
\cite{BGR}.

A function $f:D\rightarrow\mathbb{C}$ is called {\bf absolutely
continuous on lines}, written $f\in {\rm ACL}$ if for every closed
rectangular $R$ in $D$ whose sides are parallel to the coordinate
axes, $f|R$ is absolutely continuous on almost every linear segment
in $R$ which is parallel to the sides of $R$ (see, e.g., \cite{Ah},
p. 23).

Let $Q(z) : D \to I = [1, \infty]$ be an arbitrary function. A
sense-preserving homeomorphism $f : D \to \mathbb{C}$ of the class
$ACL$ is called $Q(z)-${\bf quasiconformal} ($Q(z)-$q.c.) mapping if
a.e.
\begin{equation}\label{eq3}
K_{{\mu}_{f}}(z):=\frac{1+|\mu_{f}(z)|}{1-|\mu_{f}(z)|} \leq Q(z)
\end{equation}
where ${\mu}_{f}=f_{\overline{z}}/f_{z}$ if $f_{z}\neq 0$ and
${\mu}_{f}=0$ if $f_{z}=0$. The function ${\mu}_{f}$ is called a
complex characteristic and $K_{{\mu}_{f}}$ a dilatation of the
mapping $f$.

Later ${\Bbb D}:=\left\{\nu\in{{\Bbb C}}:|\nu|\,<\,1\right\}$. Let
$\mathcal{G}$ be the qroup of all linear-fractional mappings of
${\Bbb D}$ onto itself. A set $M$ in ${\Bbb D}$ is called {\bf
invariant--convex} if all sets $g (M)$, $g \in \mathcal{G}$, are
convex, see, e.g., \cite{R}, p. 636. In particular, such sets are
convex. We say that a family of comact sets in
$M(z)\subseteq\mathbb{D}$, $ z \in \mathbb{C}$ is {\bf measurable in
the parameter} $z$, if for every closed set $M_0 \subseteq
\mathbb{C}$ the set $E_0 = \{ z \in \mathbb{C} : M(z) \subseteq M_0
\}$ is measurable by Lebesgue (cf., e.g., \cite{180}). Later we use
the following notations
\begin{equation}\label{eq4}
Q_{M}(z):=\frac{1+q_{M}(z)}{1-q_{M}(z)}\,, \,\,\,\,\,\,\,
q_{M}(z):=\max\limits_{\nu\in M(z)}|\nu|\,.
\end{equation}

Let $M(z),\,z\in\mathbb{C}$ be a family of compact sets in
$\mathbb{D}$ measurable in the parameter $z$. Let us denote by
$\mathfrak{M}_M$ the class of all measurable functions satisfying
the condition $\mu(z) \in M(z)$ a.e., and by $H_{M}^{*}$ the
collection of all regular homeomorphisms $f: \overline{\mathbb{C}}
\to \overline{\mathbb{C}}$ with the complex characteristics in
$\mathfrak{M}_M$ and the normalization $f(0)=0, \ f(1)=1, \
f(\infty)= \infty$. In the previous papers \cite{Lo} and \cite{Ram},
it was proved a series of criteria for the compactness of the
classes $H_{M}^{*}$ under the corresponding conditions on the
function $Q_{M}$, cf. also \cite{RS}, for invariant-convex
$M(z),\,z\in\mathbb{C}$. Note that the last condition implies
convexity of the set of the complex characteristics
$\mathfrak{M}_M$. As we will see later, the last circumstance
essentially simplifies the construction of variations in the classes
$H_{M}^{*}$.

One of the significant applications of compactness theorems is the
theory of the variational method. The matter is  that, in the
compact classes, it is guaranteed the existence of extremal mappings
for every continuous, in particular, nonlinear functionals. The
variational method of the research of extremal problems for
quasiconformal mappings was first applied by Belinskii~P.P., see
\cite{Be}. This method had a further development in papers of
Gutlyanskii~V.Ya., Krushkal'~S.L., Kuhnau~R., Ryazanov~V.I.,
Schiffer~M., Schober~G. and others, see, e.g.,
\cite{Gut1}--\cite{GR2}, \cite{Krush}, \cite{KrushKu1}, \cite{Rya1},
\cite{SS}, \cite{206}.

Recall that a mapping $f:X\rightarrow Y$ between metric spaces $X$
and $Y$ is called {\bf Lipschitz} if $ {\rm
dist}(f(x_{1}),\,f(x_{2}))\leq M\cdot{\rm dist}(x_{1},\,x_{2}) $ for
some $M<\infty$ and for all $x_{1},\,x_{2}\in X$ where ${\rm
dist\,}(x_{1},\,x_{2})$ denotes a distance in the metric spaces $X$
and $Y$ (see, e.g., \cite{Fe}, p. 63). The mapping $f$ is called
{\bf be-Lipschitz} if in addition $ M^{*}\cdot{\rm
dist}(x_{1},\,x_{2})\leq{\rm dist}(f(x_{1}),\,f(x_{2})) $ for some
$M^{*}>0$ and for all $x_{1},\,x_{2}\in X$.

\section{Preliminaries} Let us give necessary facts from the theory of composition operators in Sobolev's spaces.
Let $D$ be a domain in the Euclidean space $\mathbb{R}^{n}$. Recall
that the Sobolev space $L_{p}^{1}(D),\,p\geq1$, is the space of
locally integrable functions $\varphi:D\rightarrow\mathbb{R}$ with
the first partial generalized derivatives and with the seminorm
\begin{equation}\label{eq5}
\|\varphi\|_{L_{p}^{1}(D)}=\|\bigtriangledown\varphi\|_{L_{p}(D)}=
\left(\int\limits_{D}|\bigtriangledown\varphi|^{p}dm\right)^{1/p}<\infty
\end{equation}
where $m$ is the Lebesgue measure in $\mathbb{R}^{n}$,
$\bigtriangledown\varphi$ is the {\bf generalized gradient} of the
function $\varphi$,
$\bigtriangledown\varphi=\left(\frac{\partial\varphi}{\partial
x_{1}},...,\frac{\partial\varphi}{\partial x_{n}}\right)$,
$x=(x_{1},...,x_{n})$, defined by the conditions
\begin{equation}\label{eq6}
\int\limits_{D}\varphi\cdot\frac{\partial \eta}{\partial x_{i}}\,dm=
-\int\limits_{D}\frac{\partial \varphi}{\partial x_{i}}\cdot \eta\,
dm\,\,\,\,\,\,\,\,\,\,\,\,\,\,\,\,\forall\,\eta\in
C_{0}^{\infty}(D),\,i=1,\,2,...,n.
\end{equation}
As usual, here $C_{0}^{\infty}(D)$ denotes the space of all
infinitely smooth functions with a compact support in $D$.
Similarly, they say that a vector-function belongs to the Sobolev
class $L_{p}^{1}(D)$ if every its coordinate function belongs to
$L_{p}^{1}(D)$. It is known the following fact, see \cite{Uh} and
\cite{UV}.

\begin{lemma}
{}\label{l1} {\it\, Let $f$ be a homeomorphism between domains $D$
and $D^{\prime}$ in $\mathbb{R}^{n}$. Then the following statements
are equivalent:

1) the composition rule $f^{\ast}\varphi=\varphi\circ f$ generates
the bounded operator
\begin{equation}\label{eq7}
f^{\ast}:L_{p}^{1}(D^{\prime})\rightarrow L_{q}^{1}(D),\,\,\,1\leq
q\leq p<\infty\,,
\end{equation}

2) the mapping $f$ belongs to the class $W_{loc}^{1,1}(D)$ and the
function
\begin{equation}\label{eq8}
K_{p}(x,\,f):=\inf\left\{k(x):|Df|(x)\leq
k(x)|J_{f}(x)|^{\frac{1}{p}}\right\}
\end{equation}
belongs to $L_{r}(D)$ where $r$ is defined from relation
$1/r=1/q-1/p$.}
\end{lemma}

In particular, for $n=2$, $p=2$ and $q=1$, we have from here the
following statement that will be useful for us.
\begin{proposition}{}\label{p1} {\it
Let $f:\mathbb{C}\rightarrow\mathbb{C}$ be a sense-preserving
homeomorphism in the class $W_{loc}^{1,1}$ with $K_{\mu_{f}}\in
L_{loc}^{1}$. Then $g\circ f\in W_{loc}^{1,1}$ for every mapping
$g:\mathbb{C}\rightarrow\mathbb{C}$ in the class $W_{loc}^{1,2}$.}
\end{proposition}

As well-known, every quasiconformal mapping
$g:\mathbb{C}\rightarrow\mathbb{C}$ belongs to the class
$W_{loc}^{1,2}$, see, e.g., Theorem IV.1.2 in \cite{LV}. Thus, we
come to the following conclusion.

\begin{corollary}
\label{c1} {\it For every quasiconformal mapping
$g:\mathbb{C}\rightarrow\mathbb{C}$ and every sense-preserving
homeomorphism $f:\mathbb{C}\rightarrow\mathbb{C}$ of the class
$W_{loc}^{1,1}$ with $K_{\mu_{f}}\in L_{loc}^{1}$, the composition
$g\circ f$ belongs to the class $W_{loc}^{1,1}$.}
\end{corollary}

\medskip
The following statement on differentiability of the composition is
proved similarly to Theorem 5.4.6 in \cite{GR}.

\begin{lemma}
{}\label{l2} {\it Let $f$ be a homeomorphism between domains $D$ and
$D^{\prime}$ in $\mathbb{R}^{n}$, the composition operator
$f^{\ast}:L_{p}^{1}(D^{\prime})\rightarrow L_{q}^{1}(D),\,\,\,1\leq
q\leq p<\infty$, be bounded and let $f$ has $N^{-1}$--property. Then
for every function $\varphi\in L_{p}^{1}(D^{\prime})$, a.e.
\begin{equation}\label{eq9}
\frac{\partial(\varphi\circ f)}{\partial
x_{i}}(x)=\sum\limits_{k=1}^{n}\frac{\partial\varphi}{\partial
y_{k}}\left(f(x)\right)\cdot\frac{\partial f_{k}}{\partial
x_{i}}(x),\,\,\,i=1,...,n.
\end{equation}
}
\end{lemma}

Combining Lemmas \ref{l1} and \ref{l2}, similarly to IC(1) in
\cite{Ah}, we obtain.
\begin{proposition}{}\label{p2} {\it
Let $f:\mathbb{C}\rightarrow\mathbb{C}$ be a sense-preserving
regular homeomorphism with $K_{\mu_{f}}\in L_{loc}^{1}$. Then, for
every mapping $g:\mathbb{C}\rightarrow\mathbb{C}$ of the class
$W_{loc}^{1,2}$, a.e.
\begin{equation}\label{eq10}
(g\circ f)_{z}=(g_{w}\circ f)f_{z}+(g_{\overline{w}}\circ
f)\overline{f_{\overline{z}}}\,, \,\,\,\,\, (g\circ
f)_{\overline{z}}=(g_{w}\circ
f)f_{\overline{z}}+(g_{\overline{w}}\circ f)\overline{f_{{z}}}\,.
\end{equation}
}
\end{proposition}

\begin{corollary}
\label{c2} {\it In particular, formulas (\ref{eq10}) hold for
quasiconformal mappings $g:\mathbb{C}\rightarrow\mathbb{C}$.}
\end{corollary}

\section{The construction of variations} This section is devoted to
constructing variations in the classes $H_{M}^{*}$ with a method
whose idea was first proposed by Gutlyanskii V.Ya. in the paper
\cite{Gut} for analytic functions with a quasiconformal extension.
Later, this approach was applied in \cite{Rya} under constraints for
$Q_{M}$ in measure of the exponential type.

\begin{theorem}\label{t1}
{\it\, Let $M(z),\,z\in\mathbb{C}$ be an arbitrary family of convex
sets in $\mathbb{D}$. Now, let $\mu\in \mathfrak{M}_{M}$ be a
complex characteristic of a mapping $f\in H_{M}^{*}$ such that
$K_{\mu}\in L_{loc}^{1}$ and $\nu\in \mathfrak{M}_{M}$ such that the
function
\begin{equation}\label{eq11}
\varkappa=(\nu-\mu)/(1-|\mu|^{2})
\end{equation}
belongs to the open unit ball in $L^{\infty}(\mathbb{C})$. Then
there is a variation $f_{\varepsilon},\,\varepsilon\in [0,\,1/2]$ of
the mapping $f$ in the class $H_{M}^{*}$ with the complex
characteristic
\begin{equation}\label{eq12}
\mu_{\varepsilon}=\mu+\varepsilon(\nu-\mu)=(1-\varepsilon)\mu+\varepsilon\nu\,,\,\,\,\,\,\,\varepsilon\in
[0,\,1/2]
\end{equation}
such that
\begin{equation}\label{eq13}
f_{\varepsilon}(\zeta)=f(\zeta)-\frac{\varepsilon}{\pi}\int\limits_{\mathbb{C}}(\nu(z)-\mu(z))\,\varphi(f(z),\,f(\zeta))\,f_{z}^{2}\,dm_{z}+o(\varepsilon,\,\zeta)
\end{equation}
where $o(\varepsilon\,,\zeta)/\varepsilon\rightarrow0$ locally
uniform with respect to $\zeta\in \mathbb{C}$ and
\begin{equation}\label{eq14}
\varphi(w,\,w^{\prime})=\frac{1}{w-w^{\prime}}\cdot\frac{w^{\prime}}{w}\cdot\frac{w^{\prime}-1}{w-1}\,.
\end{equation}
}
\end{theorem}

{\it Proof.} Denote by $B$ a (Borel) set of all points
$z\in\mathbb{C}$ where $f$ has a total differential and
$J_{f}(z)\neq 0$. Then by definition of the class $H_{M}^{*}$ and by
the Gehring--Lehto--Menshoff theorem $|\mathbb{C}\setminus B|=0$
(see \cite{Me}, cf. Theorem III.3.1 in \cite{LV}). Moreover, by
Lemma 3.2.2 in \cite{Fe} the set $B$ can be splitted into a
countable collection of sets $B_{l}$ where $f$ is bi-Lipschitz. By
the Kirsbraun-McSchane theorem, see, e.g., Theorem 2.10.43 in
\cite{Fe}, see also \cite{Ki} and \cite{McSh}, the restrictions
$f|_{B_{l}}$ admit extensions to Lipschitz mappings of $\mathbb{C}$.
Thus, $f$ has ($N$)--property on the set $B$ and we may replace
variables in integrals, see, e.g., Theorem 3.2.5 in \cite{Fe}. Let
\begin{equation}\label{eq15}
\varkappa_{\varepsilon}
=\frac{\varepsilon\varkappa}{1-\varepsilon\varkappa\overline{\mu}}=
\varepsilon\varkappa\sum\limits_{n=0}^{\infty}\left(\varepsilon\,\varkappa\,\overline{\mu}\right)^{n},\,\,\varepsilon\in
[0,\,1]\,.
\end{equation}
Since $\|\varkappa\|_{\infty}=k<1$,
\begin{equation}\label{eq16}
\|\varkappa_{\varepsilon}\|_{\infty}\leq\frac{\varepsilon
k}{1-\varepsilon k}\leq\frac{k}{2-k}=q<1\,,\quad\ \ {\rm for}\
\varepsilon\in [0,\,1/2]\,.
\end{equation}
Now, let
\begin{equation}\label{eq17}
\gamma_{\varepsilon}(w):= \left\{
\begin{array}{rr}
\left(\varkappa_{\varepsilon}\cdot\frac{f_{z}}{\overline{f_{z}}}\right)\circ
f^{-1}(w)\,, & w\in f(B)\,,\\
0, & w\in f(\mathbb{C}\setminus B)\,.
\end{array}
\right.
\end{equation}
Re-defining, in the case of necessity, $\varkappa$ in a set of
measure zero, without loss of generality, we may assume that
$|\varkappa(z)|\leq k$ and $|\varkappa_{\varepsilon}(z)|\leq q$ for
all $z\in\mathbb{C}$ and, thus, $\gamma_{\varepsilon}(z)\leq q$ also
for all $z\in\mathbb{C}$. Moreover, since $|\mathbb{C}\setminus
B|=0$,
\begin{equation}\label{eq18}
\gamma_{\varepsilon}\circ f
=\varkappa_{\varepsilon}\cdot\frac{f_{z}}{\overline{f_{z}}}
\quad\textrm{a.e.}
\end{equation}

Consider the family of $Q$--quasiconformal
$\left(Q=(1+q)/(1-q)\right)$ mappings
$g_{\varepsilon}:\overline{\mathbb{C}}\rightarrow\overline{\mathbb{C}},\,\varepsilon\in
[0,\,1/2]$ with the complex characteristics
$\gamma_{\varepsilon},\,\varepsilon\in [0,\,1/2]$ and the
normalization $g_{\varepsilon}(0)=0$, $g_{\varepsilon}(1)=1$ and
$g_{\varepsilon}(\infty)=\infty$, see the existence theorem for
quasiconformal mappings, e.g., in the book \cite{Ah}, p. 98. By the
theorem on differentiability of $Q$--q.c. mappings in a parameter
(see \cite{Ah}, p. 105):
\begin{equation}\label{eq19}
g_{\varepsilon}(w^{\prime})=w^{\prime}-\frac{\varepsilon}{\pi}\int\limits_{f(B)}\gamma(w)
\varphi(w,\,w^{\prime})\,dm_{w}+o(\varepsilon,\,w^{\prime})
\end{equation} where $o(\varepsilon,\,w^{\prime})/\varepsilon\rightarrow 0$
as $\varepsilon\rightarrow 0$ loccally uniform with respect to
$w^{\prime}\in \mathbb{C}$ and
\begin{equation}\label{eq20}
\gamma(w)= \left\{
\begin{array}{rr}
\left(\varkappa\cdot\frac{f_{z}}{\overline{f_{z}}}\right)\circ
f^{-1}(w)\,, & w\in f(B)\,,\\
0, & w\in f(\mathbb{C}\setminus B)\,.
\end{array}
\right.
\end{equation}

Next, consider the family of mappings
$f_{\varepsilon}=g_{\varepsilon}\circ f$, $\varepsilon\in
[0,\,1/2]$. Let us show that $f_{\varepsilon}\in H_{M}^{*}$. First,
by Corollary \ref{c1}, $f_{\varepsilon}\in W_{loc}^{1,\,1}$. Then
note that the regular homeomorphism $f$ has $N^{-1}$--property by
the Ponomarev theorem, see \cite{Po}. Hence, similarly to IC(6) in
\cite{Ah}, since $J_{f}(z)\neq 0$ a.e. and $f_{z}\neq 0$ a.e. we
obtain that a.e.
\begin{equation}\label{eq21}
\mu_{g_{\varepsilon}}\circ
f=\frac{f_{z}}{\overline{f_{z}}}\cdot\frac{\mu_{f_{\varepsilon}}-\mu_{f}}{1-\overline{\mu_{f}}\cdot\mu_{f_{\varepsilon}}}\,.
\end{equation}
Here we applied the rule of differentiability of composition
(\ref{eq10}), see Corollary \ref{c2}. Solving (\ref{eq21}) with
respect to $\mu_{f_{\varepsilon}}$, we conclude that a.e.
\begin{equation}\label{eq22}
\mu_{f_{\varepsilon}}=\frac{\mu_{g_{\varepsilon}}\circ
f+\frac{f_{z}}{\overline{f_{z}}}\cdot
\mu_{f}}{\frac{f_{z}}{\overline{f_{z}}}+\overline{\mu_{f}}\cdot\mu_{g_{\varepsilon}}\circ
f}=\frac{\mu+\frac{\overline{f_{z}}}{{f_{z}}}\cdot\gamma_{\varepsilon}\circ
f}{1+\overline{\mu}\cdot\frac{\overline{f_{z}}}{{f_{z}}}\cdot\gamma_{\varepsilon}\circ
f}\,.
\end{equation}
Putting in (\ref{eq22}) the expressions from (\ref{eq15}) and
(\ref{eq18}), we have that a.e.
\begin{equation}\label{eq23}
\mu_{f_{\varepsilon}}=\frac{\mu+\varkappa_{\varepsilon}}{1+\overline{\mu}\varkappa_{\varepsilon}}=
\frac{\mu+\frac{\varepsilon\varkappa}{1-\varepsilon\varkappa\overline{\mu}}}{1+\overline{\mu}\cdot\frac{\varepsilon\varkappa}{1-\varepsilon\varkappa\overline{\mu}}}=
\mu+\varepsilon\varkappa\left(1-|\mu|^{2}\right)\,.\end{equation} By
(\ref{eq23}) and (\ref{eq11}) we obtain that
$\mu_{f_{\varepsilon}}=\mu_{\varepsilon}$ where $\mu_{\varepsilon}$
is given by (\ref{eq12}). Thus, $\mu_{f_{\varepsilon}}\in
\mathfrak{M}_{M}$, $\varepsilon\in [0,\,1/2]$ in view of convexity
of $\mathfrak{M}_{M}$.

Note that the homeomorphism $f_{\varepsilon}$ is regular for
$\varepsilon\in [0,\,1/2]$. Indeed, let us assume that
$f_{\varepsilon}$  is not regular for some $\varepsilon\in
[0,\,1/2]$. Since $|\mu_{f_{\varepsilon}}|<1$ a.e., that would be
meant that
$(f_{\varepsilon})_{z}=0=(f_{\varepsilon})_{\overline{z}}$ on a set
$E\subseteq\mathbb{C}$ of a positive measure where $f_{\varepsilon}$
is differentiable and $f$ is regular. Then similarly to IC(2) in
\cite{Ah}, we obtain that everywhere on $E$
\begin{equation}\label{eq24}
(g_{\varepsilon})_{w}\circ
f=\frac{1}{J_{f}}\left[(f_{\varepsilon})_{z}\overline{f_{z}}-
(f_{\varepsilon})_{\overline{z}}\,\overline{f_{\overline{z}}}\right]=0\,,
\end{equation}
see Proposition \ref{p2}. However, the set $\mathcal{E}:=f(E)$ has
measure zero because $g_{\varepsilon}$ is a quasiconformal mapping.
Thus, we come to the contradiction with the $N^{-1}$--property of
the mapping $f$, see \cite{Po}. Consequently, $f_{\varepsilon}\in
H_{M}^{*}$, $\varepsilon\in [0,\,1/2]$.

Finally, changing variables in (\ref{eq19}), we come to (\ref{eq13})
because $|\mathbb{C}\setminus B|=0$.

\section{Variational maximum principle} They say that a functional
$\Omega : H_{M}^{*} \to \mathbb{R}$ is {\bf differentiable by
Gateaux} if
\begin{equation}\label{eq25}
    \Omega (f_{\varepsilon}) = \Omega (f) + \varepsilon Re \int\limits_{\mathbb{C}} g\, d\varkappa + o(\varepsilon)
\end{equation}
for every variation $f_{\varepsilon} = f + \varepsilon g + o
(\varepsilon)$ in the class $H_{M}^{*}$ where $\varkappa=
\varkappa_f$ is a finite comlex Radon measure with a compact support
and  $o (\varepsilon)/\varepsilon\rightarrow 0$ as
$\varepsilon\rightarrow 0$ locally uniform in  $\mathbb{C}$ (see
\cite{206}, pp. 138--139). In other words, there is a continuous
linear functional $L (g; f)$ in the first variable such that
\begin{equation}\label{eq26}
    \Omega (f_{\varepsilon}) = \Omega (f) + \varepsilon \, Re \, L(g; f) + o
    (\varepsilon)\,.
\end{equation}

Later we assume that the function $\varphi (w, f(\zeta))$ is locally
integrable for every $f \in H_{M}^{*}$ with respect to the product
of measures $d m_w \otimes d \varkappa (\zeta)$ where $\varphi$ is
the kernel from (\ref{eq14}), $m$ is the Lebesgue measure in
$\mathbb{C}$ and that
\begin{equation}\label{eq27}
    A(w) = \frac{1}{\pi} \int\limits_{\mathbb{C}} \varphi (w, f(\zeta))\, d \varkappa (\zeta) \neq 0 \quad\textrm{for a.e.}\,\, w\in\mathbb{C}\,.
\end{equation}
Then we say that $\Omega$ is differentiable by Gateaux {\bf without
degeneration} on the class $H_{M}^{*}$.

\begin{theorem}\label{t2}
{\it Let $M(z),\,z\in\mathbb{C}$, be a family of compact convex sets
in $\mathbb{D}$ which is measurable in the parameter $z$ such that
$Q_{M}\in L_{loc}^{1}$ and let a functional $\Omega : H_{M}^{*} \to
\mathbb{R}$ is differentiable by Gateaux without degeneration. If
$\max \Omega$ in the class $H_{M}^{*}$ is attained for a mapping
$f$, then its complex characteristic satisfies the inclusion
\begin{equation}\label{eq277}
\mu(z) \in \partial M(z) \quad\textrm{for a.e.}\,\, z \in
\mathbb{C}\,.
\end{equation}}
\end{theorem}

{\it Proof.} Since $\mu \in \mathfrak{M}_M$, without loss of
generality we may assume that $\mu (z) \in M(z)$ for all $z \in
\mathbb{C}$. Let us assume that the set $$E = \{ z \in \mathbb{C} :
\mu (z) \not \in
\partial M(z) \}$$ has a positive Lebesgue measure. Let $$E_m = \{ z
\in \mathbb{C} : Q_{M} (z) \leq m \},\,\,\,\,m=1, \, 2, \, \dots,$$
$$K(z_0, r) = \{ z \in \mathbb{C} : |z - z_0 | \leq r \},\,\,\,\,\, z_0 \in
\mathbb{C}, \ r > 0\,,$$ $\chi, \ \chi_m, \ \chi_{z_0, r}$ are
characteristic functions of the sets $E$, $E_m$, $K(z_0, r)$,
respectively. Now, let $\alpha_n, \ n=1, \, 2, \, \dots,$ be an
enumeration of all rational numbers in $[0, 2\pi)$ and $\rho_n (z),
\ n=1, \, 2, \, \dots$, be distances from $\mu(z)$ till the points
of intersections of rays $\mu(z) + t e^{i \alpha_n}, \ t>0,$ with
$\partial M(z)$.

Let us show that the functions $\rho_n (z), \ n=1, \, 2, \, \dots$ ,
are measurable in the parameter $z$. Indeed, let $\Lambda_n(z) = \{
\nu \in \mathbb{C} : \nu = \mu (z) + t e^{i \alpha_n} , \quad 0 \leq
t \leq 2 \}$ be the segment of the ray passing from the point
$\mu(z)$ in the direction $e^{i \alpha_n}$ of the length 2. The
measurability of the families of the sets $\Lambda_n(z)$ in $z$
follows, e.g., from Proposition 3.1 in \cite{Ram} and general
properties of elementary operations with measurable functions (see,
e.g., \cite{180}). Consequently, the families $M_n (z) = M(z) \cap
\Lambda_n(z)$ and $\{ \eta_n(z) \} =
\partial \mathbb{D} \cap \Lambda_n (z)$ where $\partial \mathbb{D} = \{ \eta
\in \mathbb{C} : |\eta| = 1 \}$ is the unit circle are also
measurable (see Lemma 3.3 in \cite{Ram}). Thus, the functions
$\eta_n(z), \ n=1, \, 2, \, \dots$, are measurable, e.g., by the
criterion 6) in Proposition 6 in \cite{R}. By Proposition 3.1 in
\cite{Ram} the distance functions $ r_n (z) = \min\limits_{\nu \in
M_n (z)} |\nu - \eta_n (z)|$ are also measurable. It remains to note
after this that $\rho_n(z) = | \mu(z) - \eta_n(z) | - r_n (z)$.

Next, consider the functions $ \mu_n (z) = \mu (z) + \rho_n(z) e^{i
\alpha_n}$. By construction they belong to the class
$\mathfrak{M}_M$. Since the sets $M(z)$ are convex, the functions
$$\nu_n(z): = \mu(z) + \lambda(z) (\mu_n(z) -
    \mu(z))=(1-\lambda(z))\mu(z)+\lambda(z)\mu_{n}(z)$$
also belong to the class $\mathfrak{M}_M$ for an arbitrary
measurable function $\lambda (z) : \mathbb{C} \to [0, 1]$. In
particular, the class $\mathfrak{M}_M$ contains the functions
$$\nu_{z_0, r}^{m,n}(z): = \mu(z) + \lambda_m(z)\chi_{z_0, r}(z)
(\mu_n(z) - \mu(z))$$ where $$\lambda_m (z) = \frac{1-|\mu(z)|^2}
{2} \chi(z) \chi_m (z)\,.$$ Note that
$$|\mu_{n}(z)-\mu(z)|=\rho_{n}(z)\leq2q_{M}(z)$$ and
$$
  \varkappa_{z_0, r}^{m,n}(z): = \frac{\nu_{z_0, r}^{m,n} (z) - \mu(z)}{1 - |\mu(z)|^2} =
%\\ &=&
\frac{\mu_n(z) - \mu(z)}{2} \chi(z) \chi_m(z) \chi_{z_0, r} (z)
$$
belong to the closed ball of the radius $q_m: = (m-1) / (m+1) < 1$
in $L^{\infty} (\mathbb{C})$.

Since $f$ is extremal, applying the variation of Theorem \ref{t1}
with $\nu=\nu_{z_0, r}^{m,n}$, we obtain that
\begin{equation}\label{eq28}
    Re \int\limits_{\mathbb{C}} \left[ \quad \int\limits_{|z-z_0| \leq r} \varphi_{m,n} (z, \zeta)\, dm_{z} \right]
    d \varkappa (\zeta) \geq 0
\end{equation}
where $$ \varphi_{m,n} (z, \zeta) = \lambda_m (z) (\mu_n (z) -
\mu(z)) f^2_z \varphi (f(z), f(\zeta))\,.$$ Consider the functions
$$
\psi_{z_0, r}^{m,n} (w, \zeta) = \left\{
\begin{array}{rr}
\left( \varkappa_{z_0, r}^{m,n} \cdot \frac{f_z}{\overline{f_z}}
\right) \circ f^{-1} (w)\,\, \varphi (w, f(\zeta))\,, & w\in f(B)\,,\\
0, & w\in f(\mathbb{C}\setminus B)\,,
\end{array}
\right.$$ where $B$ denotes the (Borel) set of all points in
$\mathbb{C}$ where the mapping $f$ has a total differential and
$J_{f}(z)\neq0$. They are integrable with respect to the product of
the measures $d m_w \otimes d \varkappa (\zeta)$. Note that
$$J_{f^{-1}}(w)=\left[J_{f}\left(f^{-1}(w)\right)\right]^{-1}=
\left[(1-|\mu|^{2})f_{z}^{2}\right]^{-1}\left(f^{-1}(w)\right)$$ at
every point $w\in f(B)$, cf. IC(3) in \cite{Ah}. Moreover, since the
regular homeomorphism $f$ has $N^{-1}$--property, after the
replacement of variables (see Lemmas III.2.1 and III.3.2 in
\cite{LV}) we obtain that the functions $\varphi_{m,n}(z,\,\zeta)$
are also integrable with respect to the measure product $d m_z
\otimes d \varkappa (\zeta)$ and by the Lebesgue--Fubini theorem
(see, e.g., Theorem V.8.1 in \cite{51}) and (\ref{eq28}) we conclude
that
$$
\int\limits_{|z-z_0| \leq r} \left[  Re \int\limits_{\mathbb{C}}
\varphi_{m n} (z, \zeta)\,d \varkappa (\zeta)\right]
     dm_{z}  \geq 0\,.
$$

By the Lebesgue theorem on the differentiability of the indefinite
integral (see, e.g., Theorem IV(5.4) in \cite{180}) we have the
inequalities
$$\lambda_m (z) Re (\mu_n (z) - \mu(z)) \mathcal{B}(z) \geq 0 \quad\textrm{for a.e.}\,\, z \in \mathbb{C},\,\,\, m, \, n=1, \, 2, \, \dots,$$
where $\mathcal{B}(z) = \mathcal{A} (f(z)) f^2_z$ and
$\mathcal{A}(w)$ is given by (\ref{eq27}). Hence $$ \rho_n(z) Re \,
\mathcal{B}(z) e^{i \alpha_n} \geq 0,\,\,\,\,n =1, \, 2, \, \dots,
\quad\textrm{for a.e.}\,\, z \in E \cap E_m\,.$$ Since $E_m, \ m=1,
\, 2, \, \dots,$ form an exhaustion of the plane $\mathbb{C}$ in
measure, the last holds for a.e. $z \in \mathbb{C}$. On the other
hand $ \rho_n(z) > 0, \
 n =1, \, 2, \, \dots,$ on $E$ and, thus, this is equivalent to the
 inequalities
$$ Re \, \mathcal{B}(z) e^{i \alpha_n} \geq 0\,, n =1, \, 2, \,
\dots, \quad\textrm{for a.e.}\,\, z \in E\,. $$ By arbitrariness of
$\alpha_n, \ n=1, \, 2, \, \dots$, we have from here that $$ Re \,
\mathcal{B}(z) e^{i \alpha} \geq 0 \quad\forall\,\,\alpha \in [0,
2\pi) \quad\textrm{for a.e.}\,\,z \in E
 \,.$$ In particular, for $\alpha
=0$ and $\alpha=\pi$ we obtain that $ \pm Re \, \mathcal{B}(z) \geq
0, $ a.e. $ Re \, \mathcal{B}(z)= 0, $ and for $\alpha = \pi / 2$
and $\alpha= 3 \pi /2$: $ \pm \, Im \, \mathcal{B} (z) \geq 0, $
i.e. $ Im \, \mathcal{B}(z) = 0. $ Thus, $ \mathcal{B}(z)= 0$ for
a.e. $z \in E $. However, the latter is impossible because
$\mathcal{A}(w) \neq 0$ a.e., $f$ has $N^{-1}$--property and
$f_{z}\neq 0$ a.e. The obtained contradiction shows that $mes E =
0$, i.e. $\mu (z) \in
\partial M(z)$ a.e.

\section{Other necessary conditions for extremum} To formulate the
necessary conditions of the extremum we need one more notion.
Namely, let $\mu \in \mathfrak{M}_M. $ Then $\omega_{\mu} (z)$
denotes the {\bf cone of the admissible directions} (see, e.g.,
\cite{137}) for the set $M(z)$ at the point $\mu(z)$, a.e., the set
of all $\omega \in \mathbb{C}, \ \omega \neq 0,$ such that $\mu(z) +
\varepsilon \omega \in M(z)$ for all $\varepsilon \in [0,
\varepsilon_{0}]$ and some $\varepsilon_0
> 0$. Note that for strictly convex sets $M(z)$, being invariant--convex sets, the cone of admissible directions
$\omega_{\mu} (z)$ is an open cone for every $z$. Almost word for
word repeating the proof of Theorem \ref{t2}, we obtain:

\medskip

\begin{theorem}\label{t3} {\it Under the hypothesis of Theorem \ref{t2}, the extremal $f$ in the
problem on $\max \Omega$ in the class $H_{M}^{*}$ satisfies the
inequalities
\begin{equation}\label{eq29}
    Re \, \omega \mathcal{B}(z) \geq 0
\end{equation}
for a.e. $z \in \mathbb{C}$ for all $\omega$ in the cone of
admissible directions $\omega_{\mu} (z)$ where $
    \mathcal{B}(z) = \mathcal{A} (f(z))f^2_z
$ and $A(w)$ is given by (\ref{eq27}). }
\end{theorem}

\medskip

\begin{corollary}\label{c3} {\it If in addition, the boundary  is regular for a.e. $z \in
\mathbb{C}$, i.e.,  $\partial M(z)$ has a tangent at every point,
then (\ref{eq29}) is transformed to the inequality
\begin{equation}\label{eq30}
    n(z) \mathcal{B}(z) \geq 0\quad\textrm{a.e.}
\end{equation}
where $n(z)$ is the unit vector of the inner normal to $\partial
M(z)$ at the point $\mu(z)$.}
\end{corollary}

In particular, if $M(z)$ is a family of disks
\begin{equation}\label{eq31}
    M(z) = \{ \varkappa \in \mathbb{C} : |\varkappa - c(z) | \leq k(z)
    \}
\end{equation}
where the functions $c(z)$ and $k(z)$ are measurable, then by the
maximum principle, Theorem \ref{t2},  $n(z) = (c(z) - \mu(z))/k(z)$
and the relation (\ref{eq30}) is equivalent to the equality
$$\frac{c(z) - \mu(z)}{k(z)} =
\frac{\overline{\mathcal{B}(z)}}{|\mathcal{B}(z)|}\quad\textrm{a.e.}\,,
$$ i.e., $$\mu(z)= c(z) - k(z)
\frac{\overline{\mathcal{B}(z)}}{|\mathcal{B}(z)|}\,. $$ Thus, we
have:

\vskip7pt

\begin{corollary}\label{c4} {\it Let $M(z), \ z \in \mathbb{C}$, be the family of the disks (\ref{eq31}),
$k+|c|\in L_{loc}^{1}$, and the functional $\Omega : H_{M}^{*} \to
\mathbb{R}$ is differentiable by Gateaux without degeneration. Тhen
the extremal of the problem on $\max \Omega$ in the class
$H_{M}^{*}$ satisfies the equality
\begin{equation}\label{eq32}
    f_{\overline{z}}= c(z) f_z - k(z) \frac{\overline{\mathcal{A}(f(z))}}{|\mathcal{A}(f(z))|}\, \overline{f_z}\,.
\end{equation}
In particular, if $c(z) =0$ we obtain the equality
\begin{equation}\label{eq33}
    f_{\overline{z}}= - k(z)\, \frac{\overline{\mathcal{A}(f(z))}}{|\mathcal{A}(f(z))|}\,\, \overline{f_z}\,.
\end{equation}
}
\end{corollary}

\medskip
\noindent {\it Tatyana Lomako, Vladimir Ryazanov} \\
Institute of Applied Mathematics and Mechanics,\\
National Academy of Sciences of Ukraine \\
74 Roze Luxemburg str., 83114 Donetsk, Ukraine \\
%Phone: +38 -- (062) -- 3110145, Fax: +38 -- (062) -- 3110285 \\
E-mail: tlomako@yandex.ru, vlryazanov1@rambler.ru

\end{document}